\documentclass[12pt]{amsart}
\usepackage{amsmath,amssymb,latexsym,amsthm,graphicx,citesort}
\newtheorem{theorem}{Theorem}
\newtheorem{lemma}[theorem]{Lemma}

\newtheorem{corollary}[theorem]{Corollary}

\def\qed{$\Box$}

\def\RR{\mathbb{R}}
\def\ZZ{\mathbb{Z}}

\begin{document}
\title{Range conditions for a spherical mean transform}
\author{Mark Agranovsky}
\address{{\bf Permanent address:} Mathematics Department\\ Bar Ilan University\\
Ramat Gan 52900
Israel\\ {\bf Spring 2009:} Mathematics Department\\
Texas A\&M University\\
College Station, TX 77843-3368, USA}
\email{agranovs@macs.biu.ac.il}
\author{David Finch}
\address{Mathematics Department\\Oregon State University\\
Corvallis, Oregon 97331-4605, USA}
\email{finch@MATH.ORST.EDU}
\author{Peter Kuchment}
\address{Mathematics Department\\
Texas A\&M University\\
College Station, TX 77843-3368, USA}
\email{kuchment@math.tamu.edu}
\date{}
\maketitle
\begin{abstract}
The paper is devoted to the range description of the Radon type transform that averages a function over all spheres centered on a given sphere. Such transforms arise naturally in thermoacoustic tomography, a novel method of medical imaging. Range descriptions have recently been obtained for such transforms, and consisted of smoothness and support conditions, moment conditions, and some additional orthogonality conditions of spectral nature. It has been noticed that in odd dimensions, surprisingly, the moment conditions are superfluous and can be eliminated. It is shown in this text that in fact the same happens in any dimension.
\end{abstract}

\section{Introduction}

This papers deals with the integral transform that associates to a function $f \in C(\RR^n)$ its spherical means
$$
g(x,t)=\int_{|\theta|=1} f(x+t\theta) d\theta, \quad x\in\RR^n, t\geq 0,
$$
where $d\theta$ is the normalized area measure on the unit sphere
$$
\{\theta \in \RR^n:|\theta|=1\}.
$$

This integral geometric transform arises in a variety of pure and applied mathematics topics from PDEs, medical imaging, and geophysics \cite{John,CH,AKQ,AQ,Beylkin,LP1,LP2,FPR}, which we will not survey here.

Experience of integral geometry and its applications (e.g., \cite{Helg_Radon,GGG,GGG1,GelfVil,Leon_Radon,Natt_new}) shows importance of the range descriptions for integral geometric Radon type operators. For the spherical mean transform $f\mapsto g$ introduced above, the range is provided by the famous Asgeirsson's theorem \cite{Asgeirsson,John,CH,Helg_Radon}, which says that the resulting functions $g(x,t)$
are exactly the solutions of the following initial value problem in $\RR^n$ for the Darboux-Poisson-Euler equation:
\begin{equation}\label{E:Darboux}
\begin{cases}
\frac{\partial ^2 g}{\partial t^2}+\frac{n-1}{t}\frac{\partial
g}{\partial t}=\Delta g\\
g(x,0)=f(x)\\
\frac{\partial g}{\partial t}(x,0)=0.
\end{cases}
\end{equation}

Note, however, that the spherical mean operator $f\mapsto g$ increases the number of variables from $n$ to $n+1$. It is thus natural to reduce the choice of spheres to an $n$-dimensional family. This can be done in various ways. One is to allow only a discrete set of values of the radius $t$. This kind of a restricted spherical mean transforms has been intensively studied for many years in harmonic analysis (see \cite{Volch,Zalc} and references therein) and will not be discussed here. Our target transform is different, the one where the set of centers $x$ is restricted to an $(n-1)$-dimensional variety $S\subset \RR^n$. This restricted type of the operator was studied, due to its relevance to approximation theory, in \cite{LP1,LP2,AQ}. We in fact restrict this transform even further, assuming $S$ to be the unit sphere
$$
S:=\{x \in \RR^n:|x|=1\}.
$$
Moreover, it will be assumed that the original function $f$ is supported in the unit ball
$$
B:=\{x \in \RR^n:|x|\leq 1\}.
$$
Thus, our article is devoted to the study of the following {\it restricted spherical mean operator} $R_S$:
$$
g(x,t):=R_S f(x,t):=\int_{|\theta|=1} f(x+t\theta) d\theta, \quad (x,t) \in S \times \RR^+.
$$
The transform just introduced has been studied intensively in recent years (see surveys \cite{AKK,PatchSch,FR3,KuKu}), since it represents the data collected in, and thus is the major component of the newly developing method of medical imaging called thermoacoustic tomography
\cite{Wang_book,MXW_review}.

Our interest here is in obtaining the range descriptions (i.e., analogs of Paley-Wiener theorem) for the operator $R_S$ on the space $C^\infty_0(B)$.

It is instructional to look at the example of the standard Radon transform $R$ that produces integrals $g(\theta,s)$ of $f(x)$ over $(n-1)$-dimensional planes $x\times\theta=s$. The well known \cite{Helg_Radon,Leon_Radon,GGG1,Natt_new} range description for $R$ contains (besides some natural smoothness requirements) three types of conditions:
\begin{itemize}
\item the natural support condition on $g=Rf$: planes not touching the support of $f$ produce zero integrals,
\item the so called moment conditions: $M_k(\theta):=\int s^k g(\theta,s)ds$ is the restriction to the unit sphere of a homogeneous polynomial of degree $k$, and
\item the evenness condition $g(\theta,s)=g(-\theta,-s)$ that reflect the fact that the planes given by the equations $x\cdot\theta=s$ and $x\cdot(-\theta)=-s$ coincide.
\end{itemize}
On the face of it, the moment conditions constitute the only non-trivial part. However, this is where the standard Radon transform is misleading. Indeed, for various its generalizations the situation changes. For instance, for some weighted (so called attenuated) Radon transforms arising in emission tomography, the analogs of moment conditions were found rather fast, while it took many years for an analog of the ``trivial'' evenness condition to be found \cite{KucLvin,Nov}. A similar thing has happened with the transform $R_S$: moment type conditions were explicitly described in \cite{Patch}, although they were implicitly present already in \cite{LP1,LP2,AQ}. It was harder to find the complete set of conditions, i.e. to find analogs of the evenness condition. This was done in $2d$ in \cite{AmbKuc_range}, then in odd dimensions in \cite{FR}, and finally various range descriptions that apply to any dimension were obtained in \cite{AKQ}. In most cases, range descriptions were in a sense parallel to the Radon case: they consisted of a natural support condition for $g=R_S f$, the moment conditions, and some ``orthogonality'' condition that thus looked like an analog of evenness in the Radon case. However, the following peculiar thing was observed in odd dimensions in \cite{FR}: the moment conditions could be dropped completely from the range description. In other words, the orthogonality conditions were apparently not just analogs of evenness, but also contained the moment conditions as well. This absolutely non-transparent fact was later confirmed in the odd dimensions for the different range descriptions obtained in \cite{AKQ}. The initial impression of the authors was that the situation might depend on the validity of the Huygens' principle for the wave equation, which only holds in odd dimensions \cite{CH}. However, the goal of this text is to show that {\bf the moment conditions can be dropped from the range descriptions in any, not only odd, dimension}.

The next Section \ref{S:result} contains the formulation of the result, which is then proven in Section \ref{S:proof}. The paper ends with the Remarks and Acknowledgments sections and the bibliography.

\section{Formulation of the result}\label{S:result}

Let us first introduce some notations.

We use the standard notation $J_p(\lambda)$ for Bessel functions of the first kind. We will also need to use what is sometimes called {\em normalized} or {\em spherical} Bessel functions \cite{Magnus}:
\begin{equation}\label{E:modif_Bessel}
j_p (\lambda)=\frac{2^p\Gamma (p+1) J_p (\lambda)}{\lambda^p}.
\end{equation}
For $p=\frac{n-2}{2}$, according to Poisson's representations of
Bessel functions \cite{John,Magnus}, this is the spherical
average of a plane wave in $\RR^n$.

Consider for each $m\in\ZZ^+$ an orthonormal basis  $Y^m_l, 1 \leq l \leq d(m)$ in the space of all spherical harmonics of degree $m$ in $\RR^n$.

Let $g(x,t)$ be a function on $S\times \RR^+$. This function can be expanded into the Fourier series with respect to spherical harmonics as follows:
\begin{equation}\label{E:sph_series}
    g(x,t)=\sum\limits_{m,l}g_{m,l}(t)Y^m_l(x),
\end{equation}
where
\begin{equation} \label{E:gml}
 g_{m,l}(t)=\int_{\theta \in S} g(\theta,t) Y^m_l(\theta) d\theta.
\end{equation}

We will also need the Fourier-Bessel transform of $g_{m,l}$:
\begin{equation*}
\widehat g_{m,l}(\lambda)=\int_0^{\infty} g_{m,l}(t) j_{n/2-1}(\lambda t) t^{n-1}dt.
\end{equation*}
The main result of this paper is the following:
\begin{theorem}\label{T:Main}
    The following three statements for a function $g\in C^\infty_0 (S\times [0,2])$ are equivalent:
    \begin{enumerate}
    \item The function $g$ is
    representable as $R_Sf$ for some $f\in C^\infty_0(B)$.

    \item For any eigenvalue $-\lambda^2$ of the Dirichlet Laplacian in $B$
    and the corresponding eigenfunction $\psi_\lambda (x)$, the following orthogonality condition is satisfied:
    \begin{equation}\label{E:orthog}
\int\limits_{S\times [0,2]} g(x,t)\partial_\nu \psi_\lambda (x)j_{n/2-1}(\lambda t)t^{n-1}dxdt=0.
    \end{equation}
    Here $\partial_\nu$ is the exterior normal derivative at the
    boundary of $B$.

    \item   For any integer $m$ and $1\leq l \leq d(m)$, function $\widehat{g}_{m,l}(\lambda)$ vanishes at all
    non-zero zeros of the Bessel function $J_{m+n/2-1}(\lambda)$ (equivalently, at all zeros of $j_{m+n/2-1}(\lambda)$).
     \end{enumerate}
\end{theorem}

Not all the claims of this Theorem are new. Indeed, implication $(1) \Rightarrow (2),(3)$, as well as equivalence $(2)\Leftrightarrow (3)$ were established in \cite[Theorem 10]{AKQ}. In fact, it was also shown there that each of the statements (2) and (3) implies (1), if being supplemented by additional
\begin{description}\label{C:Moment}
  \item[Moment conditions:]

  for any non-negative integer $k$, the function $M_k$ on $S$ defined as
\begin{equation}\label{E:Moment}
M_k(x):=\int\limits_0^2 t^{2k+n-1} g(x,t)dt, \quad x\in S
\end{equation}
has an extension to  $\RR^n$ as a polynomial of degree at most $2k$.
\end{description}

Therefore, due to \cite[Theorem 10]{AKQ}, the proof of Theorem \ref{T:Main} will be concluded, if we prove the following surprising result:
\begin{theorem}\label{T:moments}
For a function $g\in C^\infty_0 (S\times [0,2])$, any of the equivalent conditions (2) and (3) implies the moment conditions (\ref{E:Moment}).
\end{theorem}

\section{Proof of Theorem \ref{T:moments}}\label{S:proof}

Since the orthogonality condition 3 of Theorem \ref{T:Main} is formulated in the Fourier domain, we start the proof with reformulating the moment conditions in a similar way, as it was done in \cite{AKQ}. The next lemma is proven in \cite[Lemma 9]{AKQ}. For completeness, we provide its brief proof here.

\begin{lemma}\label{L:order} \cite[Lemma 9]{AKQ})
The following conditions on a function $g\in C^\infty_0(S\times[0,2])$ are equivalent:
\begin{enumerate}
\item The moment conditions (\ref{E:Moment}) for $g(x,t)$.

\item For any $\ m \geq 1$ and  $1\leq l \leq d(m)$, function $\widehat g_{m,l}(\lambda)$ vanishes at $\lambda=0$ to the order at least $m$.
\end{enumerate}
\end{lemma}
\proof
Expanding $g(\theta,t)$ into an the
orthonormal basis $Y^m_l$ of spherical harmonics on $S$, one obtains
\begin{equation}\label{E:momenta_harmonics}
   M_k(\theta)= \int\limits_0^\infty t^{2k+n-1}g(\theta,t)dt=\sum\limits_{l,m} Y^m_l(\theta)\int\limits_0^\infty t^{2k+n-1}g_{m,l}(t)dt.
\end{equation}
A spherical harmonic $Y^m_l$ of degree $m$ can be extended to a polynomial
of degree $d$ if and only  if $d\geq m$ and $d-m$ is even. Thus,
the moment conditions require that
$$\int_0^\infty t^{2k+n-1}g_{m,l}(t)dt=0
$$
for $m>2k$.

On the other hand, using the well known series expansion $j_p(\lambda)=\sum_{k=0}^\infty C_k\lambda^{2k}$ for Bessel functions, one arrives to the formula
\begin{equation}\label{E:g_hat_expansion}
\widehat g_{m,l} (\lambda)=\sum\limits_{k=0}^\infty C_k\lambda^{2k}
\int\limits_0^\infty t^{2k+n-1}\int\limits_{\theta \in
S}Y^m_l(\theta)g(\theta,t)dS(\theta) dt .
\end{equation}
This shows that the moment conditions are equivalent to the
vanishing of all terms in this series with $2k<m$.
\qed

We now return to the proof of Theorem \ref{T:moments}. We assume that $g\in C^\infty_0(S\times [0,2])$ satisfies conditions 3 of Theorem \ref{T:Main}, i.e. that for any $m\geq 1, 1\leq l \leq d(m)$, function $\widehat{g}_{m,l}(\lambda)$ vanishes at all zeros of $j_{m+n/2-1}(\lambda)$. We need to show that then $\lambda=0$ is a zero of order at least $m$ of $\widehat{g}_{m,l}(\lambda)$.

In fact, we will achieve this working with each individual
function $\widehat{g}_{m,l}(\lambda)$. To simplify notations, from now on, when considering a functions $g_{m,l}$ or $\widehat g_{m,l}$, we will drop in the notations the irrelevant index $l$, which does not enter any considerations. Thus, $g_m$ will denote any of the functions $g_{m,l}$ with $1\leq l \leq d(m)$.

So, let $g$ satisfy the conditions above. Then $\widehat g_m(\lambda)$ is an entire function of the following Paley-Wiener class:
\begin{equation}\label{E:gPW}
|\widehat g_m(\lambda)|\leq C_N (1+|\lambda|)^{-N}e^{2|Im\, \lambda|} \mbox{ for any $N>0$.}
\end{equation}
We also know that it vanishes at all zeroes of the entire function $j_{m+n/2-1}(\lambda)$, which happen to be simple. Thus, the ratio $H(\lambda)$ of these two functions is also entire:
\begin{equation}\label{E:divis}
\widehat g_m(\lambda)=j_{m+n/2-1}(\lambda) H(\lambda).
\end{equation}
We will now determine the Paley-Wiener class to which $H$ belongs. This is done using the following simple (and probably well known) estimate from below for the Bessel
function of the first kind $J_\nu$:
\begin{lemma}\label{L:bessel}\cite{AmbKuc_range}
On the entire complex plane, except for a disk $S_0$ centered at
the origin and a countable number of disks $S_k$ of radii $\pi/6$
centered at points $\pi(k+\frac{2\nu+3}{4})$, one has
\begin{equation}
\label{E:estimate} |J_\nu (z)|\ge\frac{C
e^{|Im\,z|}}{\sqrt{|z|}},\quad C>0.
\end{equation}
\end {lemma}
For completeness, we present here the proof of this lemma, borrowed from \cite{AmbKuc_range}.

\textbf{Proof of Lemma \ref{L:bessel}:} We split the complex plane into three parts by
a circle $S_0$ of a radius $R$ (to be chosen later) centered at
the origin and a planar strip $\{z=x+iy|\,|y|<a\}$. Here part I consists
of points $z$ such that $|z|\geq R$ and
$|Im\,z|\geq a$; in part II, $|z|\geq R$
and $|Im\,z|< a$; finally, part III is the interior of $S_0$, i.e. contains points $z$ such that $|z|<R$.
It is sufficient to prove the estimate (\ref{E:estimate})
in the first two parts: outside and inside the strip. In fact, it suffices to consider only the right
half plane $Re\,z\geq 0$ intersected with parts I and II.

The Bessel function $J_\nu(z)$ has the following
known asymptotic representation in the sector $|\arg
z|\le\pi-\delta$ (e.g., \cite[formula (4.8.5)]{Andrews} or
\cite[formula (5.11.6)]{Lebedev}):
\begin{equation}\label{E:asymptotic}
\begin{array}{c}
J_\nu(z)=\sqrt{\frac{2}{\pi z}} \cos(z-\frac{\pi \nu}{2}
-\frac{\pi}{4})(1+O(|z|^{-2}))\\
 -\sqrt{\frac{2}{\pi z}} \sin(z-\frac{\pi
\nu}{2}-\frac{\pi}{4})\left( \frac{4\nu^2-1}{8z}+O(|z|^{-3})\right)
\end{array}
\end{equation}

Let us estimate it in the part I,
where $|Im \, z|>a$ and $|z|>R$ for sufficiently large $a$
and $R$ (and, as we have agreed, $Re\,z \geq 0$). There one has
$\displaystyle{\frac{\sin z}{z} = \cos z\; ( O(|z|^{-1}))}$,
and thus (\ref{E:asymptotic}) implies
$$
J_\nu(z)=\sqrt{\frac{2}{\pi z}} \cos(z-\frac{\pi \nu}{2}
-\frac{\pi}{4})(1+O(|z|^{-1})).
$$
For sufficiently large $a, R$ this leads to
the needed estimate
\begin{equation}
\label{E:part1} |J_\nu(z)|\ge\frac{C e^{|Im\, z|}}{\sqrt {|z|}}
\end{equation}

In the part II of the plane (right half of the strip), due to
boundedness of $\sin z$, one has
$$
J_\nu(z)=\sqrt{\frac{2}{\pi z}}\left[\cos(z-\frac{\pi \nu}{2}
-\frac{\pi}{4})(1+O(|z|^{-2}))+ O(|z|^{-1})\right].
$$
Let us consider the system of non-intersecting circles $S_k$ of radii $\frac{\pi}{6}$ centered
at $z_k=\frac{\pi}{2}+k\pi+\frac{\pi \nu}{2}+\frac{\pi}{4}$. Outside these circles,
$|\cos(z-\frac{\pi \nu}{2} -\frac{\pi}{4})|\ge C$ and thus
$$
|J_\nu(z)|\ge
\frac{C}{\sqrt{|z|}}(1+O(|z|^{-1})),
$$
which implies, inside of the strip and outside the circles $S_k$, the estimate
\begin{equation}
\label{E:part2} |J_\nu(z)|\ge\frac{C e^{|Im\, z|}}{\sqrt{|z|}}
\end{equation}
for $|z|>R$ with a suitably chosen sufficiently large
$R$. This proves the statement of the lemma. \qed

\begin{corollary}\label{C:PW}
    Function  $H(\lambda)$ in (\ref{E:divis}) is even, entire, and satisfies the following Paley-Wiener estimate
    \begin{equation*}
        | H(\lambda)|\leq C_N(1+|\lambda|)^{-N}e^{|Im\, \lambda|} \mbox{for any $N>0$.}
    \end{equation*}
\end{corollary}
\proof Since $H=\widehat g_m /j_{m+n/2-1}$, this follows from (\ref{E:gPW}) and (\ref{E:estimate}). \qed

The well known Paley-Wiener theorem for the Fourier-Bessel transform (see, e.g. \cite[Lemma 1]{AKQ} and references in that paper) shows that then the following representation is possible:
\begin{equation}
H(\lambda)=\int_0^{\infty} h(t) j_{(n-2)/2} (t) t^{n-1}dt,
\end{equation}
where $h(t)\in C^{\infty}_0(\RR)$ is an even function on the real line supported in the interval $[-1,1]$.

We introduce now the function $\widetilde h$ on $\RR^n$ as
$$
\widetilde h(y)=h(|y|), \mbox{ where }|y|=\sqrt{\sum y_j^2}.
$$
Then $\widetilde h$ is clearly a $C^{\infty}$-function on $\RR^n$ supported in the unit ball $B(0,1)$ centered at the origin.

We would like to lift similarly the Bessel function $j_{m+n/2-1}$ to $\RR^n$. This produces
the function
$$
E_m(v):=j_{m+n/2-1}(|v|)
$$

\begin{lemma}\label{L:E_m} Function $E_m(v)$ is, up to a constant factor, the $n$-dimensional Fourier transform of the function $(1-|y|^2)_+^{m-1}, y \in \RR^n $, i.e.
\begin{equation}\label{E:E_m}
E_m(v) = \mbox{const} \int_{\RR^n} e^{-i v \cdot y} (1-|y|^2)_+^{m-1} dy.
\end{equation}
Here $u_+:=max(u,0)$.
\end{lemma}

\proof This identity follows from the Sonine integral formula \cite[formula 4.11.11]{Andrews}, but we will indicate the proof, to make the presentation more self-contained.

The integral on the right hand side is clearly invariant under rotations, so
we may suppose that $v = \vert v\vert e_1$. We set $y =(s,y')$ with $ y'
\in \RR^{n-1}.$ Applying Fubini's theorem, the integral on the right can
then be written
\[ \int_{-1}^{1} e^{-i\vert v\vert s}\int_{s=\vert v\vert} (1 - s^2 -\vert
y'\vert^2)_+^{m-1} \,dy' \, ds.\] The substitution $y' = \sqrt{1-s^2}
\sigma$ reduces this to \[ C_{m,n}\int_{-1}^{1} e^{-i\vert v\vert s}
(1-s^2)^{m + n/2-1 -1/2} \,ds,\] with $C_{m,n}= \frac{1}{2}\vert S^{n-2}\vert
B(\frac{n-1}{2},\frac{m}{2}),$ where $\vert S^{n-2}\vert $ is the area of the $(n-2)$-sphere,
and $B(\cdot,\cdot)$ is Euler's Beta function.
The integral is a constant multiple of a well known
representation of $j_{m+n/2-1}.$ \qed

\begin{corollary} Function $ g_m(|v|)$ on $\RR^n$ is representable as the convolution
\begin{equation}\label{E:conv}
g_m (|v|)= w_m(|v|) * \widetilde h,
\end{equation}
where $w_m(|v|)=(1-|v|^2)_+^{m-1}$ and $\widetilde h \in C_0^{\infty}(\RR^n)$ is supported in the unit ball $B(0,1)$.
\end{corollary}
Indeed, this follows from (\ref{E:divis}) and Lemma \ref{L:E_m}.

Let us now pay attention to the behavior of the function $g_m$ at the point $t=0$. Namely, since  $g_m(t)\in C^\infty[0,2]$, it vanishes at $t=0$ with all its derivatives. This behavior is inherited by the smooth function $g_m(|v|)$ in $\RR^n$: it vanishes at the origin with all its derivatives. In particular, if $D=(\partial/\partial v_1,... \partial/\partial v_n)$ and $\alpha=(\alpha_1, ... ,\alpha_n)$ is a multi-index, we get the  identity
\begin{equation}\label{E:D}
0= D^\alpha g_m(0)=   w_m*\left( D^{\alpha} \widetilde
h\right) \vert_{v=0}=
\int\limits_{\RR^n} (1-|y|^2)_+^{m-1}\left(D^{\alpha}\widetilde h(y)\right)dy.
\end{equation}
Since $\widetilde h$ is smooth and supported in the unit ball, the right
hand side is equal to \[ \int\limits_{\RR^n} (1-\vert y\vert^2)^{m-1}
\left(D^{\alpha}\widetilde h(y)\right)\,dy.\]
Integrating by parts, we now conclude that the following equalities hold for any $\alpha$:
\begin{equation}\label{E:Polyn}
\int\limits_{\RR^n}\left(D^\alpha(1-|y|^2)^{m-1}\right)\widetilde h(y)dy=0.
\end{equation}

We can now establish the claim central to the proof of Theorem \ref{T:moments}:
\begin{lemma}\label{L:Laplace} For any $k \leq m-1$, the following moment relation holds:
\begin{equation}
\int\limits_{\RR^n} |y|^{2k} \widetilde h(y) dy=\int\limits_{|y|<1} |y|^{2k} \widetilde h(y) dy=0.
\end{equation}
\end{lemma}
\proof Indeed, expanding $(1-|y|^2)^{m-1}$ according to the binomial formula and using (\ref{E:Polyn}), leads in particular to
$$
\sum\limits_{k=0}^{m} \binom{m}{k}
(-1)^k \int\limits_{|y|<1} \Delta^j (|y|^{2k}) \widetilde h(y)dy=0,
$$
for any integer $j\geq 0$, where $\Delta^j$ is the $j$-th power of the $n$-dimensional Laplace operator. Using the well known and easily verifiable relation $\Delta |y|^{2k}=c_k|y|^{2k-2}$ with $c_k\neq 0$, and running the value $j$ from $j=0$ to $j=m-1$, one obtains a uniquely solvable triangular homogeneous linear system of equations with respect to the moments
$\int\limits_{|y|<1} (|y|^{2k}) \widetilde h(y)dy$ for $0\leq k \leq m-1$. Thus, these moments are all equal to zero. Due to the evenness of $h$, all odd order moments also vanish.
\qed

The vanishing of the moments shown in the last Lemma implies that the Fourier transform $h$ of $\widehat h$, and hence the corresponding one-dimensional function $H$, has zero at the origin of order at least $m$. Due to (\ref{E:divis}), this feature is inherited by $\widehat g_m$, which finishes the proof of Theorem \ref{E:Moment}, and thus also of Theorem \ref{T:Main}.

\section{Remarks}

\begin{enumerate}
\item The implication $(1) \mbox{ in Theorem \ref{T:Main} } \Rightarrow \mbox{ moment conditions}$ was implicitly present in \cite{AQ,LP1,LP2} and explicitly formulated as range conditions in \cite{Patch}.

\item If $S$ is not a sphere, most probably the moment conditions are independent of the orthogonality ones.

\item We have shown that in fact $\widehat g_{m,l}$ has at the origin a zero of the order $2m$, not just mere $m$. As it was mentioned in a similar situation in \cite{AKQ}, this feature reflects the fact that $S$ is a sphere and most probably will not hold for more general surfaces $S$.

\item The range conditions provided in \cite{AKQ} had also another incarnation, besides the versions (2) and (3) of Theorem \ref{T:Main}. Indeed, it was shown in \cite[Theorem 3]{AKQ} that the following condition is equivalent to (2) and (3):\\
{\em The backword in time solution of the Darboux equation (\ref{E:Darboux}) with the boundary data provided by the function $g$ and zero initial conditions for a sufficiently large moment $T$ of time does not develop a singularity when $t\to 0$}. (see details in \cite{AKQ}).

This condition can be substituted instead of (2) or (3) in Theorem \ref{T:Main}.
\end{enumerate}

\section*{Acknowledgments}
The work of the first author was partially supported by the ISF (Israel Science Foundation) Grant 688/08 and by the Texas A\&M University. The third author was partially supported by the NSF grant DMS 0604778 and by the KAUST grant KUS-CI-016-04. The authors express their gratitude to NSF, Texas A\&M University, and KAUST for the support. The authors also thank L.~Nguyen and Rakesh for useful discussions and comments.

\end{document}